\documentclass[10pt]{book}
\usepackage{array}
\usepackage{color}
\usepackage{subfigure}
\usepackage{amsthm,amsmath}
\usepackage{remreset}
\usepackage{graphicx}
\usepackage{enumerate}
\usepackage[all]{xy}
\usepackage{hyperref}

\makeatletter
\@removefromreset{table}{chapter}
\@removefromreset{figure}{chapter}
\makeatother

\newtheorem{thm}{Theorem}
\newtheorem{defi}{Definition}

\newtheorem{property}{Property}
\newtheorem{conclusion}{Corollary}
\newtheorem{example}{Example}

\begin{document}

\def\letas{\mathrel{\mathop{=}\limits^{\triangle}}}
\def\ind{\begin{picture}(9,8)
         \put(0,0){\line(1,0){9}}
         \put(3,0){\line(0,1){8}}
         \put(6,0){\line(0,1){8}}
         \end{picture}
        }
\def\nind{\begin{picture}(9,8)
         \put(0,0){\line(1,0){9}}
         \put(3,0){\line(0,1){8}}
         \put(6,0){\line(0,1){8}}
         \put(1,0){{\it /}}
         \end{picture}
    }

\renewcommand{\baselinestretch}{1.2}

\markright{
\hbox{\footnotesize\rm }\hfill
}

\markboth{\hfill{\footnotesize\rm Zhichao Jiang, Peng Ding, and Zhi Geng} \hfill}
{\hfill {\footnotesize\rm Transitivity of association signs} \hfill}

\renewcommand{\thefootnote}{}
$\ $\par


\fontsize{10.95}{14pt plus.8pt minus .6pt}\selectfont
\vspace{0.8pc}

\vspace{5mm}
\centerline{\large\bf
Qualitative evaluation of associations by the transitivity of }

\centerline{\large\bf the association signs }

\vspace{.4cm}
\centerline{Zhichao Jiang$^1$, Peng Ding$^2$ and Zhi Geng$^1$}
\vspace{.4cm}
\centerline{\it $^1$Peking University and $^2$Harvard University}
\vspace{.55cm}
\fontsize{9}{11.5pt plus.8pt minus .6pt}\selectfont


\begin{quotation}
\noindent {\it Abstract:}
We say that the  signs of association measures 
among three variables $\{X, Y,  Z\}$ are transitive
if a positive association measure between the variable $X$ and  the intermediate variable $Y$
and further a positive association measure between $Y$ and the endpoint variable $Z$ imply a positive association measure between $X$ and $Z$.  We introduce four association measures with  different stringencies, and discuss conditions for the transitivity of the signs of these association measures. When the variables follow exponential family distributions, the conditions become simpler and more interpretable. Applying our results to two data sets from
an observational study and a randomized experiment, we demonstrate that the results can help us to draw conclusions about the signs of the association measures between $X$ and $Z$ based only on two separate studies  about $\{ X, Y\}$ and $\{ Y, Z\}$.
\par

\vspace{9pt}
\par
\end{quotation}\par


\fontsize{10.95}{14pt plus.8pt minus .6pt}\selectfont

\vspace{5mm}
\noindent {\bf 1.  Introduction }

Reasoning by transitivity is commonly, at least implicitly, applied
to statistical results from different studies.
For association measures, however, the transitivity may not be guaranteed without conditions.
For example, suppose an epidemiologic study found
that irregular heart beat
had a positive association with sudden death,
and we also know from a clinical trial
that a certain drug  significantly corrects irregular heart beat.
Might we conclude from these two statistical results that the drug can reduce the rate of sudden death?
What conditions are required  for reasoning from the known statistical results?
Transitivity and its required conditions are important
for reasoning from the known results of association measures  cumulated by statistical inferences.
For example, in a meta-analysis, some associations between variable pairs may be obtained from published papers, but the original data may not be available.
What conditions are required to qualitatively evaluate
an association between another pair of variables based on these known associations?
There are a few statistical approaches for such qualitative reasonings.
VanderWeele and Robins (2010)
propose the signed directed acyclic graph (DAG) to qualitatively reason the sign of association or  causal measure
between two variables
by the signs of directed edges on the paths between the two variables.
Their approach requires both a whole DAG over all variables
and the signs of all edges on each path between the two variables.
VanderWeele and Tan (2012) propose an approach for propagation of bounds within the DAG framework, which requires a known DAG with edge-specific bounds.
Rubin (2004) discusses how to combine the results from two randomized trials of anthrax vaccine  on human volunteers and macaques, where the outcome of interest, survival when challenged with lethal doses of anthrax, is available only from macaques.
Pearl and Bareinboim (2012) discuss
the transportability of causal effects across different studies
based on the DAG framework.
Prentice (1989) proposes a criterion for surrogate endpoints in clinical trials
so that a null effect of treatment on a surrogate implies
a null effect of treatment on the endpoint.
Chen et al. (2007) and Ju and Geng (2010)
discuss the transitivity of causal effects
between variable pairs of treatment, surrogate and endpoint.

In this paper, we focus on the transitivity of the signs of association measures. We introduce
 four association measures with different
stringency levels: density, cumulative distribution, expectation and correlation levels
(Cox and Wermuth, 2003; Whittaker, 1990).
We say that the signs of association measures are transitive if a non-negative
(or positive) association between $X$ and $Y$ and a non-negative (or
positive) association between $Y$ and $Z$ imply a non-negative (or positive)
association between $X$ and $Z$. We discuss the transitivity of the signs of
association measures, and present conditions and assumptions (or prior
knowledge) required for the transitivity.
We show that a more stringent association measure has stronger transitivity.
We focus on the transitivity of the signs of association measures
among three variables, and these results can be easily extended to cases
with more variables.
We discuss conditions for the transitivity of these association measures separately
for two cases with and without the conditional independence of
$X$ and $Z$ given $Y$.
Conditional independence is  one condition of Prentice's criterion for
evaluating surrogacy of $Y$
for the endpoint $Z$ (Prentice, 1989).
The conditions for transitivity proposed in this paper allow for qualitative assessment of
the association between two variables $X$ and $Z$,
 by the data sampled from the marginal
distributions of $(X,Y)$ and $(Y,Z)$ or from the conditional distribution of $(X,Z)$ given $Y$.

The remainder of the paper is organized as follows.
Section 2 presents the definitions of four association measures and
discusses their stringencies.
In Section 3, we consider the transitivity of the signs of these association measures
under the conditional independence of $X$ and $Z$ given $Y$, and give  results for an exponential family distribution.
We generalize the results about transitivity without conditional independence in Section 4. We apply our theoretical results to two data sets in Section 5;
we conclude with a discussion in Section 6, and give
all proofs of the theorems in the Appendix.

\vspace{5mm}
\noindent {\bf 2. Association measures and their stringencies}

In this section, we introduce four commonly-used association measures and show their relative stringencies for depicting the associations.

\begin{defi}
The four association measures between $X$ and $Y$ are:
\begin{enumerate}[(1)]
\item  density association:
$\partial^2 \ln f(x, y) / \partial x \partial y$ (Whittaker, 1990);
\item  distribution association:
$\partial F(y|x) / \partial x$ (Cox and Wermuth, 2003);
\item  expectation association:
$\partial E(Y|x) / \partial x$;
\item  correlation coefficient: $r(X, Y)$.
\end{enumerate}
\end{defi}

For these measures, $X$ and $Y$ may be continuous, discrete or mixed random variables.
For a discrete variable, we can replace the partial differentiation
by the difference between two adjacent levels.
For instance, when $X$ and $Y$ are both binary variables,
the density association is the log odds ratio
$$\ln \frac{P(X=1, Y=1)P(X=0, Y=0)}{P(X=1, Y=0)P(X=0, Y=1)},$$
 and the distribution association and  expectation association are both equal to the risk difference $P(Y=1 | X=1)-P(Y=1 | X=0)$. Moreover, if we are concerned only about the signs of the association measures, a non-negative expectation association implies that the risk ratio is greater than or equal to one, i.e. $P(Y=1 | X=1)/P(Y=1 | X=0) \geq 1$.

The density association $\partial^2 \ln f(x, y) / \partial x \partial y$ depicts a local dependence between $X$ and $Y$ around the point $(x,y)$.
An important property of the density association is
\begin{eqnarray*}
\frac{\partial^2 \ln f(x, y)}{\partial x \partial y}
= \frac{\partial^2 \ln f(x| y)}{\partial x \partial y}
= \frac{\partial^2 \ln f(y|x)}{\partial x \partial y} ,
\end{eqnarray*}
which can be identified by sampling
conditionally on $X$ or $Y$, such as a case-control sampling
or a follow-up study.
The distribution association $\partial F(y|x) / \partial x$ depicts the dependence of a global $Y\leq y$ on a local $X=x$,
and $\partial F(y|x) / \partial x \leq 0$
means that $Y$ given $X$ is stochastically increasing in $X$ (Cox and Wermuth, 2003).
The expectation association $\partial E(Y|x) / \partial x$ depicts the overall dependence of $Y$ on $X$,
and correlation coefficient $r(X,Y)$ depicts a linear association.
We say that $X$ and $Y$ are non-negatively associated
with respect to a measure, say $\partial^2 \ln f(x, y) / \partial x \partial y$,
if $\partial^2 \ln f(x, y) / \partial x \partial y \geq 0$ for all $x$ and $y$,
and we say that $X$ and $Y$ are positively associated with respect to a measure if further strict inequality holds for some $x$ or $y$.
Let $A \ind B$ denote the independence between $A$ and $B$, and let $A \ind B |C$ denote the conditional independence of $A$ and $B$ given $C$.
Generally, a non-negative association measure at a more stringent level
implies a non-negative association measure at a less stringent level,
as summarized in the following properties:

\begin{property}

The implication relationship (Xie et al., 2008):
\begin{eqnarray*}
& & \frac{\partial^2 \ln f(x, y)}{\partial x \partial y}\geq0, \forall x, y  \Longrightarrow \frac{\partial F(y|x)}{\partial x}\leq0, \forall x, y \Longrightarrow  \frac{\partial E(Y|x)}{\partial x}\geq0, \forall x   \\ & & \Longrightarrow  r(X, Y) \geq 0;
\end{eqnarray*}
\end{property}
\begin{property}
The equivalent relationship among null association measures:
$$
   \frac{\partial^2 \ln f(x, y)}{\partial x \partial y}=0, \forall x, y
   \Longleftrightarrow \frac{\partial F(y|x)}{\partial x}=0, \forall x, y
   \Longleftrightarrow X \ind Y;
$$
\end{property}

\begin{property}
For a bivariate normal vector $(X,Y)$, the equivalent relationship:
\begin{eqnarray*}
\frac{\partial^2 \ln f(x, y)}{\partial x \partial y}\geq0, \forall x, y & \Longleftrightarrow & \frac{\partial F(y|x)}{\partial x}\leq0, \forall x, y \Longleftrightarrow  \frac{\partial E(Y|x)}{\partial x}\geq0, \forall x  \\ & \Longleftrightarrow & r(X, Y) \geq 0;
\end{eqnarray*}
\end{property}

\begin{property}
For a binary $Y$, the equivalent relationship:
$$
\frac{\partial^2 \ln f(x, y)}{\partial x \partial y}\geq0, \forall x, y \Longleftrightarrow  \frac{\partial F(y|x)}{\partial x}\leq0,  \forall x, y \Longleftrightarrow  \frac{\partial E(Y|x)}{\partial x}\geq0,  \forall x;
$$
\end{property}

\begin{property}
For a binary $X$, the equivalent relationship:
$$
\frac{\partial E(Y|x)}{\partial x}\geq0,  \forall x \Longleftrightarrow r(X, Y)\geq 0.
$$
\end{property}
Note that all the  relationships above
are also true for the strict inequalities and equalities, that is,
(``$\geq$'', ``$\leq$'') in the above expressions
can be changed to (``$>$'', ``$<$'') and (``$=$'', ``$=$'').

By the implication relationship in Property 1,
 the density association is the most stringent, and the correlation coefficient is the least stringent.
For the case of two normal variables or two binary variables, all these association measures have the same signs (non-negative, null or
positive).

\vspace{5mm} \noindent {\bf 3. Transitivity of association signs with conditional independence}

Now we consider three variables $\{ X, Y, Z\}$,
and discuss the conditions for the transitivity of
association signs among them.
Prentice (1989) uses conditional independence
as a criterion for validating a surrogate $Y$, when evaluating
the effect of treatment $X$ on the endpoint $Z$.
We call $X\ind Z|Y$ the conditional independence assumption, where $Y$ breaks the dependence between $X$ and $Z$.
 Note that we can have both $X \ind Z | Y$ and $X \ind Z$ even when both pairs $(Y, Z)$ and $(X, Y)$ are associated, and Table \ref{tab:Birch} is an example of this given in Birch (1963). This is nowadays called a violation of weak transitivity or of singleton transitivity (Ln\v{e}ni\v{c}ka and Mat\'{u}\v{s}, 2007), in which case we may not infer  the transitivity 
 of associations only under the conditional independence assumption.
 It may be because the associations within different
level sets have different signs.
For example, the association between $X$ and $Y$ within $\{Y=0, Y=1\}$
has a different sign from that between $X$ and $Y$ within $\{Y=1, Y=2\}$. 
In this section, we discuss the transitivity of the association signs under the conditional independence assumption,
and we extend them to the case without the conditional independence assumption in Section 4.

\begin{table}
\begin{center}
\caption{Example of violation of weak transitivity}
\label{tab:Birch}
\begin{tabular}{cccccccccccc} \hline
 & \multicolumn{2}{c}{$Y=0$} && \multicolumn{2}{c}{$Y=1$}& &\multicolumn{2}{c}{$Y=2$} & & \multicolumn{2}{c}{all $Y$} \\
\cline{2-3} \cline{5-6}  \cline{8-9}  \cline{11-12}
      & $Z=1$& $Z=0$& & $Z=1$& $Z=0$& & $Z=1$& $Z=0$& & $Z=1$& $Z=0$ \\ \hline
$X=1$        &  4& 2& & 2&   1& & 1& 4& & 7 & 7 \\
$X=0$        &  2& 1 && 4& 2 &  & 1& 4& & 7  & 7\\ \hline
all $X$ & 6&3& & 6&3&& 2&8&& 14&14 \\ \hline
\end{tabular}
\end{center}
\end{table}

\vspace{5mm} \noindent {\bf 3.1 Transitivity without distributional assumptions}

We shall show, under the conditional independence assumption, that
the more stringent association measures,
the density and distribution associations, are transitive.
 But the other two less stringent association measures,
the expectation association and correlation, are not transitive
without additional conditions or assistance from the more
stringent measures.
The first theorem is about the transitivity of the density association.

\begin{thm} \label{th1}
(Density association)
Under the assumption $X \ind Z|Y$, if
$$
\mbox{(1) }
\frac{\partial^2 \ln f(x, y)}{\partial x \partial y} \geq0, \forall x, y
\hspace{1cm} \mbox{ and } \hspace{1cm}
\mbox{(2) }
\frac{\partial^2 \ln f(y, z)}{\partial y\partial z} \geq0, \forall y, z,
$$
then $\partial^2 \ln f(x,z) / \partial x\partial z \geq 0, \forall x, z.$
\end{thm}

According to Theorem \ref{th1},
we can conclude a non-negative (positive) density association between $X$ and $Z$
from the prior knowledge
of a non-negative (positive) density association 
between $X$ and $Y$ and a non-negative (positive)
density association between $Y$ and $Z$.
Notice that the transitivity of the signs is also applicable to
non-positive association measures
if we define $X'=-X$ or $Y'=-Y$.

Although Theorem \ref{th1} gives the result for non-negative signs, it is much more important to have the result stated directly for only positive or negative signs of associations.  Even if the conditions in Theorem \ref{th1} are satisfied, we may have $\partial^2 \ln f(x,z) / \partial x\partial z = 0$ for all $x$ and  $z$, i.e., $X \ind Z$. If strict inequalities hold for conditions (1) and (2) in Theorem \ref{th1} for any $(x,y,z)$ in
a set of nonzero measure,
then we have $\partial^2 \ln f(x,z) / \partial x\partial z > 0$
in this set, 
that is, the density association measure between $X$ and $Z$ is positive. In fact, if all the variables are discrete, we will alway obtain the strict inequality unless $X\ind Y$ or $Y \ind Z$.

Similarly, we have the following transitivity of the distribution association.

\begin{thm} \label{th2}
(Distribution association)
Under the assumption $X \ind Z|Y$, if
$$
\mbox{(1) } \frac{\partial F(y|x)}{\partial x}\leq0, \forall x, y
\hspace{1cm} \mbox{ and } \hspace{1cm}
    \mbox{(2) } \frac{\partial F(z|y)}{\partial y}\leq0, \forall y, z,
$$
then $\partial F(z|x) / \partial x \leq 0, \forall x, z$.
\end{thm}

The expectation associations themselves are not transitive,
and they require assistance from more stringent measures, as shown in the following theorem.

\begin{thm} \label{th3}
(Expectation association)
Under the assumption $X \ind Z|Y$, if
$$
\mbox{(1) } \frac{\partial F(y|x)}{\partial x}\leq0, \forall x, y
\hspace{1cm} \mbox{ and } \hspace{1cm}
\mbox{(2) } \frac{\partial E(Z|y)}{\partial y}\geq0, \forall y,
$$
then $\partial E(Z|x) / \partial x \geq 0, \forall x$.
\end{thm}

Similar to Theorem \ref{th1},
for Theorems \ref{th2} and \ref{th3},
if there exists
a set of nonzero measure
in which strict inequalities hold in (1) and (2),
then the corresponding association measures
between $X$ and $Z$ are positive.

Theorems \ref{th1} to \ref{th3} are also useful  for cases with more than three variables. If we have  another variable vector $V$, then the conditions in these Theorems should be conditional on $V$, and we can obtain the association signs of $X$ and $Z$ conditional on $V$. This is  useful for models with a covariate vector $V$.

  Notice that in Theorems \ref{th1} to \ref{th3}, if we have either $X \ind Y$ or $Y \ind Z$, we will always obtain $X \ind Z$.
Condition (1) of Theorem \ref{th3}
is a distribution association but not an expectation association
of $Y$ on $X$. Below we give a numerical example to illustrate
that the expectation association of $Y$ on $X$ cannot replace condition (1) of Theorem \ref{th3}.

\begin{example}
We generate data under conditional independence: $X\sim$ Bernoulli$(1/2), \\
\varepsilon\sim $ Bernoulli$(p), Y = X+2\varepsilon (1 - X), Z = I(Y = 2)$, where $p< 1/2$ and $I(\cdot)$ is the indicator function.
We have that $E(Y|X=1) - E(Y|X=0)=1- 2p\geq 0$, $E(Z|Y=2)-E(Z|Y=1)=1-0\geq 0$ and $E(Z|Y=1)- E(Z|Y=0) = 0$,
but we calculate that $E(Z|x=1)-E(Z|x=0)=0-p \leq 0$.
\end{example}

However, for a linear model of $Z$ given $Y$,
the expectation association measures
are transitive, and further the transitivity can be
represented by an equation of expectation associations as follows.

\begin{conclusion}
\label{co1}
Under the assumption $X\ind Z|Y$, if $E(Z|y)=\alpha+\beta y$,
then
$\partial E(Z|x) / \partial x = \beta \partial E(Y|x) / \partial x$.
\end{conclusion}
 In Example 1, we see that $Z$ does not follow a linear model given $Y$, and thus we cannot infer the transitivity of association signs.
Using the implication relationship of Property 1 and Theorems \ref{th1} to \ref{th3}, we summarize the transitivity of association signs
in Table \ref{tab1}.
We show the transitivity of non-negative association measures
between $X$ and $Y$ and between $Y$ and $Z$ to a non-negative association measure
between $X$ and $Z$.
We see from Table \ref{tab1}
that a non-negative association measure between $X$ and $Z$
requires the same or more stringent non-negative association measures
between $X$ and $Y$ and between $Y$ and $Z$.
Notice that Table \ref{tab1} is not symmetric,
and the expectation association of $Y$ on $X$
does not have any implication for the sign of an association measure between $X$ and $Z$ unless $Z$ has a linear model, as shown in the last line of Table \ref{tab1}.
In the following example, a non-negative expectation association of $Y$ on $X$
and even the most stringent non-negative density association between $Y$ and $Z$ do not imply a non-negative expectation association of $Z$ on $X$.

\begin{table}
\begin{center}
\caption{Transitivity of association signs under conditional independence} \label{tab1}
\begin{tabular}{c|ccccc} \hline
Association  &  \multicolumn{5}{c}{Association between $Y$ and $Z$} \\ \cline{2-6}
between $X$ and $Y$
& $\frac{\partial^2 \ln f(y, z)}{\partial y \partial z} \geq0, \forall y, z$ & $\Rightarrow$
& $\frac{\partial F(z|y)}{\partial y}\leq0, \forall y, z$ & $\Rightarrow$
& $\frac{\partial E(Z|y)}{\partial y}\geq0, \forall y$ \\ \hline
$\frac{\partial^2 \ln f(x, y)}{\partial x \partial y} \geq0, \forall x, y$ &
$\frac{\partial^2 \ln f(x, z)}{\partial x \partial z} \geq0, \forall x, z$ & &
$\frac{\partial F(z|x)}{\partial x}\leq0, \forall x, z$ & &
$\frac{\partial E(Z|x)}{\partial x}\geq0, \forall x$  \\
$\Downarrow$ & \\
$\frac{\partial F(y|x)}{\partial x}\leq0, \forall x, y$ &
$\frac{\partial F(z|x)}{\partial x}\leq0, \forall x, z$ & &
$\frac{\partial F(z|x)}{\partial x}\leq0, \forall x, z$ & &
$\frac{\partial E(Z|x)}{\partial x}\geq0, \forall x$  \\
$\Downarrow$ &  \\
$\frac{\partial E(Y|x)}{\partial x}\geq0, \forall x$ &
\multicolumn{5}{c}{Under $E(Z|y) = \alpha + \beta y$,
$\frac{\partial E(Z|x)}{\partial x}
= \beta \frac{\partial E(Y|x)}{\partial x}\geq 0, \forall x$.}  \\ \hline
\end{tabular}
\end{center}
\end{table}

\begin{example}
Assume $ X \ind Z|Y$ with the distributions $P(y|x)$ and $P(z|y)$ given in Table \ref{tab2}.
Then we have
$E(Y|X=1) - E(Y|X=0) = 0.2$ and
$$\ln \frac{P(Y=y, Z=0)P(Y=y+1,Z=1)}{P(Y=y,Z=1)P(Y=y+1,Z=0)} \geq 0,$$
for $y=0$ and 1, but $E(Z|X=1) - E(Z|X=0)=-0.32$.
\end{example}

\begin{table}[htbp]
    \caption{Distributions $P(y|x)$ and $P(z|y)$ for Example 2} \label{tab2}
    \centering
    \subtable[]{
    \begin{tabular}{c|ccc}\hline
    & $Y=0$ & $Y=1$ & $Y=2$ \\ \hline
    $X=0$ & 0.6 & 0 & 0.4 \\
    $X=1$ & 0 & 1 & 0 \\ \hline
    \end{tabular}
    }
    \qquad
    \subtable[]{
    \begin{tabular}{c|ccc}\hline
    & $Y=0$ & $Y=1$ & $Y=2$  \\ \hline
    $Z=0$ & 0.9 & 0.9 & 0.1 \\
    $Z=1$ & 0.1 & 0.1 & 0.9 \\  \hline
    \end{tabular}
    }
\end{table}

The transitivity of the correlation coefficient signs
is weaker than that of expectation association measure signs.
Under the conditional independence assumption,
we cannot obtain a non-negative correlation between $X$ and $Z$
from a non-negative correlation of $X$ and $Y$ ($Y$ and $Z$)
and another non-negative association measure between $Y$ and $Z$ ($X$ and $Y$).
This is illustrated by the following example.

\begin{example}
Assume $X \ind Z|Y$ with
the distributions $P(y|x)$ and $P(z|y)$ given in Table \ref{tab3}.
Then we have $cov(Y,Z)=0.0017>0$ and
$$\ln \frac{P(Y=y, X=0)P(Y=y+1,X=1)}{P(Y=y,X=1)P(Y=y+1,X=0)} \geq 0$$
for $y=0$ and 1, but $E(Z|X=1) - E(Z|X=0)=-0.005$.
\end{example}

\begin{table}[htbp]
    \caption{Distributions $P(y|x)$ and $P(z|y)$ for Example 3} \label{tab3}
    \centering
    \subtable[]{
    \begin{tabular}{c|ccc}\hline
    & $Y=0$ & $Y=1$ & $Y=2$ \\ \hline
    $X=0$ & 0.1 & 0.1 & 0.8 \\
    $X=1$ & 0.05 & 0.05 & 0.9 \\ \hline
    \end{tabular}
    }
    \qquad
    \subtable[]{
    \begin{tabular}{c|ccc}\hline
    & $Y=0$ & $Y=1$ & $Y=2$  \\ \hline
    $Z=0$ & 0.3 & 0 & 0.2 \\
    $Z=1$ & 0.7 & 1 & 0.8 \\  \hline
    \end{tabular}
    }
\end{table}

\vspace{5mm} \noindent {\bf 3.2 Transitivity in the exponential family}

In the previous subsection,
we exhibited sufficient conditions for
non-negative and positive association measures between $X$ and $Z$,
but such conditions may not be necessary.
In other words, it is possible that $X$ and $Z$ are non-negatively or positively associated, but $X$ and $Y$ (or $Y$ and $Z$)
are negatively associated.
We shall show the equivalent relationship between the sign
of association measure between $X$ and $Y$ and that between $X$ and $Z$,
under the assumptions that $Y$ follows an exponential family distribution and the association between $Y$ and $Z$ is non-negative.
Below we review the definition of the exponential family.

\begin{defi}
  We say that $Y$ given $X$ follows an exponential family distribution
  if its density (or its probability mass function for a discrete $Y$) has the form
$$
f(y|x; \theta, \phi) = \exp \left\{  \frac{ y\theta_{x} - b(\theta_{x})}{ a(\phi) }+c(y,\phi)  \right\},
$$
where $\theta_x$ is a function of $x$ and $a(\phi)>0$.
\end{defi}

In particular, a binomial or normal $Y$
is from the  exponential family.
For the exponential family,
we show an equivalent relationship
among the signs of the density, distribution and expectation associations:

\begin{thm} \label{co2}
If $Y$ given $X$ follows an exponential family distribution, then
$$
\frac{\partial^2 \ln f(x,y)}{\partial x \partial y}\geq0,\forall x,y \Longleftrightarrow \frac{\partial F(y|x)}{\partial x}\leq0 ,\forall x,y
   \Longleftrightarrow \frac{\partial E(Y|x)}{\partial x}\geq0,\forall x.
$$
The equivalent relationship is also true if (``$\geq$'', ``$\leq$'')
in the inequalities is changed to (``$>$'', ``$<$'') and (``$=$'', ``$=$'').
\end{thm}

Suppose below that we have prior knowledge of the sign of association
between $Y$ and $Z$, and we discuss the equivalent relationships
between the signs of associations between $X$ and $Y$ and between $X$ and $Z$.

\begin{conclusion} \label{co3}
Assume $X \ind Z|Y$
and that $Y$ given $X$ follows an exponential family distribution.
\begin{enumerate}[(1)]
\item
     If $\partial E(Z|y) / \partial y \geq 0, \forall y$
     and strict inequality holds for a nonzero measure set, then
$$
\frac{\partial E(Y|x)}{\partial x}\geq 0, \forall x \Longleftrightarrow
\frac{\partial E(Z|x)}{\partial x} \geq 0, \forall x;
$$

\item   If $\partial F(z|y) / \partial y \leq 0, \forall y,z$
     and strict inequality holds for a nonzero measure set, then
$$
\frac{\partial F(y|x)}{\partial x}\leq0, \forall x, y \Longleftrightarrow
\frac{\partial F(z|x)}{\partial x}\leq0, \forall x, z;
$$

\item  If $\partial^2 \ln f(y, z) / \partial y \partial z \geq 0, \forall y,z$  and strict inequality holds for a nonzero measure set, then
$$
\frac{\partial^2 \ln f(x,y)}{\partial x \partial y}\geq0, \forall x, y \Longleftrightarrow
\frac{\partial^2 \ln f(x,z)}{\partial x \partial z}\geq0, \forall x, z.
$$

\end{enumerate}
\end{conclusion}

These equivalent relationships above are also true if the inequalities (``$\geq$'', ``$\leq$'') above are changed to
strict inequalities (``$>$'', ``$<$'') and equalities (``$=$'', ``$=$'').

In practice, when direct measurement of  $Z$ is too time-consuming or costly, we may try to measure a surrogate endpoint $Y$ instead
and use the sign of the association between $X$ and $Y$ to predict the
sign of the association between $X$ and $Z$.

The results above can be extended to transitivity of  causal measures. 
If $X$ is randomized or is conditionally independent of
all potential outcomes given some covariates,
then the associations between $X$ and $Y$ and between $X$ and $Z$
are also the causal effects of $X$ on $Y$ and $X$ on $Z$, respectively.
If we have prior knowledge that 
the conditional independence assumption holds and the  association between $Y$ and $Z$ is non-negative, then according to Corollary \ref{co3}, we may conclude
 that the sign of the treatment effect of $X$ on $Z$
is the same as the sign of the treatment effect of $X$ on  $Y$.

\vspace{5mm} \noindent {\bf 4. Transitivity of association signs without conditional independence}

In many real applications, the conditional independence assumption $X \ind Z|Y$ may not hold.
In this section, we remove the conditional independence assumption,
and discuss conditions for the transitivity of these association signs.

 \vspace{5mm} \noindent {\bf 4.1 Motivating examples}

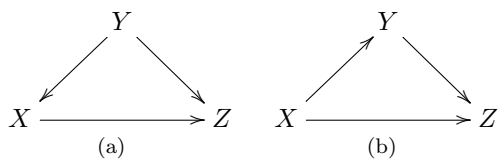
\begin{figure}
\centering

\subfigure[]{ \xymatrix{
   & Y \ar[dl] \ar[dr]&\\
 X  \ar[rr]&  & Z }}
 {\         }
\subfigure[]{\xymatrix{
   & Y \ar[dr]&\\
 X \ar[ru] \ar[rr]&  & Z }}
 \caption{Two  DAGs for $(X,Y,Z)$}
 \label{fig:dag}
\end{figure}

 In this section, we consider two cases.
In the first case as shown in Figure \ref{fig:dag}(a), 
$Y$ is a confounder between $X$ and $Z$.
In the second case as shown in Figure \ref{fig:dag}(b), 
$X$ has a direct path to $Z$ and 
an indirect path to $Z$ through an intermediate variable $Y$.

  For Figure \ref{fig:dag}(a), our results can help us to find whether the Yule--Simpson Paradox occurs, since we can evaluate whether the conditional association sign is the same as the marginal association sign between $X$ and $Z$. Appleton et al. (1996) give an example of the Yule--Simpson Paradox, and we show the condensed form in Table \ref{tab:simpson}. As shown in Table \ref{tab:simpson}, the association between $X$ and $Z$ is positive conditional on $Y$ but negative marginally, thus the  Yule--Simpson Paradox occurs.
 
\begin{table}
\centering
\caption{Numbers of women smokers and nonsmokers in different age groups.}
\label{tab:simpson}
\begin{tabular}[ht]{ccccccccccccccc}\hline
$Y$: age group & \multicolumn{2}{c}{18-34} && \multicolumn{2}{c}{35-54} && \multicolumn{2}{c}{55-64} && \multicolumn{2}{c}{65+} &&\multicolumn{2}{c}{all ages} \\
\cline{2-3} \cline{5-6}  \cline{8-9}  \cline{11-12}  \cline{14-15}
$Z/X$: smoker &yes & no &&yes & no &&yes & no &&yes & no& &yes & no \\ \hline
dead &5 &6 &&41 &19 &&51 &40& &42 &105 &&139 &230 \\
alive & 174 &213 &&198&180&&64&81&&7&28&&443 &502 \\ \hline
odds ratio & \multicolumn{2}{c}{1.02}& & \multicolumn{2}{c}{1.96}&& \multicolumn{2}{c}{1.61}&& \multicolumn{2}{c}{1.02}&& \multicolumn{2}{c}{0.68} \\ \hline
\end{tabular}
\end{table}

 For Figure \ref{fig:dag}(b), our results are useful  to draw conclusions about
 the effect of $X$ on $Z$ based on some assumptions about the pairs $(X,Y)$ and
 $(Y,Z)$. These are generalizations of the results in Section 3 which do not
 allow a directed arrow from $X$ to $Z$. As shown in Table \ref{tab:trans}, we
 have that $E(Z=1|Y=1)-E(Z=1|Y=1)=0.08$ and $E(Y=1|X=1)-E(Y=1|X=1)=0.2$, but
 $E(Z=1|X=1)-E(Z=1|X=1)=-0.08$. Thus without the conditional independence
 assumption, even if the pairs $(X,Y)$ and $(Y,Z)$ are both positively
 associated, $X$ and $Z$ may be negatively associated.

\begin{table}
\centering
\caption{Distribution $P(x, y, z)$ for violation of transitivity}
\label{tab:trans}
\begin{tabular}[ht]{ccccccccccc}\hline
 & \multicolumn{2}{c}{$Y=1$} & & \multicolumn{2}{c}{$Y=0$} \\
 \cline{2-3} \cline{5-6}
&$X=1$ & $X=0$ && $X=1$ & $X=0$ \\ \hline
$Z=1$ &  0.15 & 0.12 & & 0.08 & 0.15\\
$Z=0$ & 0.15 & 0.08 && 0.12 &0.15\\ \hline
\end{tabular}
\end{table}

\vspace{5mm} \noindent 
{\bf 4.2 Transitivity without distributional assumptions}

 Analogous to Theorems \ref{th1} to \ref{th3}, we give conditions for the transitivity of
density, distribution and expectation association signs as follows.

\begin{thm}
(Density association) \label{thm:6}
Assume $\partial^2 \ln f(x, z|y) / \partial x\partial z \geq 0, \forall x, y, z$. If
\begin{eqnarray*}
& \mbox{(1)} & \frac{\partial^2 \ln f(x, y)}{\partial x \partial y} \geq 0, \ \ \  \forall x, y, \\
& \mbox{(2)} & \frac{\partial^2 \ln f(y, z|x)}{\partial y\partial z}\geq 0, \ \ \  \forall x, y, z, \mbox{ and } \\
& \mbox{(3)}& \frac{\partial^2 \ln f(z|x, y)}{\partial x \partial y}\geq 0, \ \ \
\forall x, y, z,
\end{eqnarray*}
then we have $\partial^2 \ln f(x, z) / \partial x\partial z \geq 0$, $\forall x, z$.
\end{thm}

Instead of the conditional independence $X \ind Z|Y$, we require a prior knowledge about
the sign of association between $X$ and $Z$  given  $Y$.

If we have prior knowledge that the association between $X$ and $Z$ given $Y$
is non-positive,
we can replace $Z$ by $Z'=-Z$ in the assumptions.
In Theorem \ref{thm:6}, condition (1) means that the density association between $X$ and $Y$ is non-negative, and condition (2) means that the density association between $Y$ and $Z$ conditional on $X$ is non-negative. Condition (3) is quite different and can be interpreted as a non-negative interaction of $(X,Y)$ on $Z.$

\begin{thm} (Distribution association) \label{thm:7}
Assume $\partial F(z|y, x) / \partial x \leq 0, \forall x ,y, z.$ If
$$
\mbox{(1) } \frac{\partial F(y|x)}{\partial x} \leq 0, \ \  \forall x, y,
\hspace{1cm} \mbox{ and } \hspace{1cm}
\mbox{(2) } \frac{\partial F(z|y, x)}{\partial y} \leq 0, \ \ \forall x, y, z,
$$
then $\partial F(z|x) / \partial x \leq 0$, $\forall x, z$.
\end{thm}

\begin{thm}(Expectation association) \label{thm:8}
Assume $\partial E(Z|y, x) / \partial x \geq 0, \forall x,y$. If
$$
\mbox{(1) }
\frac{\partial F(y|x)}{\partial x}\leq 0, \ \  \forall x, y,
\hspace{1cm} \mbox{ and }
\hspace{1cm} \mbox{(2) }
\frac{\partial E(Z|y, x)}{\partial y}\geq 0, \  \ \forall x, y,
$$
then $\partial E(Z|x) / \partial x \geq 0$, $\forall x$.
\end{thm}

Comparing conditions (2) in Theorems \ref{thm:6} to \ref{thm:8} with those in Theorems \ref{th1}
to \ref{th3},  we see that
without the conditional independence assumption, the association between $Y$ and $Z$ is required to be conditional on $X$.

 There are some implications from Theorems \ref{thm:6} to \ref{thm:8}. First, we do not need the joint distribution $f(x,y,z)$ to judge whether the Yule--Simpson Paradox occurs:
\begin{conclusion}\label{co:simpson}
All the conditions in Theorem \ref{thm:6} to \ref{thm:8} can be evaluated by $f(x,z | y)$.
\end{conclusion}
Corollary \ref{co:simpson} implies that we can assess the sign of the marginal association
measure between $X$ and $Z$ by
the conditional distribution $f(z, x|y)$,  and the marginal distribution of $Y$ is not required.

Second,
by relaxing the conditional independence assumption to the assumption
of non-negative association measure between $X$ and $Z$  given  $Y$,
the marginal associations between $Y$ and $Z$
are replaced by the conditional associations given $X$. However,
this condition can be weakened if
the models of $Z$ and $Y$ are linear.

\begin{conclusion} \label{co4}
Assume that
$E(Z|x,y)=\beta_{0}+\beta_{1}x+\beta_{2}y$, $E(Y|x)=\beta_{3}+\beta_{4}x, \beta_1 \geq 0$, and $\beta_4\geq 0$. If
$$
\mbox{(1) }
\beta_2 \geq 0,
\hspace{1cm} \mbox{ or } \hspace{1cm}
\mbox{(2) }
\frac{\partial E(Z|y)}{\partial y} \geq 0, \ \ \forall y,
$$
then we have $\partial E(Z|x) / \partial x = \beta_1 + \beta_2 \beta_4 \geq 0$ for all $x$.
\end{conclusion}

If we have $\beta_2 \geq 0$, i.e., $Y$ and $Z$ are non-negatively associated conditional on $X$, then Corollary  \ref{co4} is the same as  Cochran (1938). But even if we have only the non-negative association between $Y$ and $Z$ marginally, we can still conclude that $X$ and $Z$ are non-negatively associated based on the linear models in Corollary \ref{co4}. Thus, we do not need to observe $X$ to judge condition (2) in Corollary \ref{co4}. Provided $X$ and $Z$ are positively associated conditional on $Y$, if we have two populations, one with $(X,Y)$ observed and the other with $(Y,Z)$ observed, we can draw the conclusion about the association sign between $X$ and $Z$.

Third, when there is a fully randomized intervention on $Y$, we will have $X \ind Y$, and  Theorems \ref{thm:6} and \ref{thm:8} can be further simplified.

\begin{conclusion}\label{co:ind}
Suppose $X \ind Y$.
\begin{enumerate}[(1)]
\item If $X$ or $Z$ is binary, then $\partial^2 \ln f(x,z| y)/ \partial x\partial z \geq 0, \forall x,y,z$ implies $\partial^2 \ln f(x,z)/ \partial x\partial z \geq 0, \forall x,z$.

\item If $\partial F(z|y,x)/\partial x \leq 0, \forall x,y $, then $\partial F(z|x)/\partial x \leq 0, \forall x$;

\item If $\partial E(Z|y,x)/\partial x \leq 0, \forall x,y $, then $\partial E(Z|x)/\partial x \leq 0, \forall x$.
 \end{enumerate}
\end{conclusion}

Fourth, Theorems \ref{thm:6} to \ref{thm:8} can also be useful  for the cases with more than three variables in the same way as Theorems \ref{th1} to \ref{th3}. In addition,  there is another way to use  Theorems \ref{thm:6} to \ref{thm:8}  in these cases. We can combine $X$ and $V$  into a vector
that plays the same role as the original $X$. 
It is not applicable in Theorems \ref{th1} to \ref{th3} 
since $(X,V) \ind Z | Y$ may not hold.

\vspace{5mm} \noindent {\bf 4.3 Transitivity in the exponential family}

 Combining Theorems \ref{thm:6} to \ref{thm:8} with the results for the exponential  family in Theorem \ref{co2}, we can strengthen our conclusions as follows:

\begin{thm}\label{thm:exp}
Assume that $Z$ given
$X$ follows an exponential family distribution, and that $\partial E(Z|y, x) / \partial x \geq 0, \forall x, y$. If
$$
\mbox{(1) }
\frac{\partial F(y|x)}{\partial x}\leq 0, \ \  \forall x, y,
\hspace{1cm} \mbox{ and }
\hspace{1cm} \mbox{(2) }
\frac{\partial E(Z|y, x)}{\partial y}\geq 0, \  \ \forall x, y,
$$
then $\partial^2 \ln f(x, z) / \partial x\partial z \geq 0$, $\forall x, z$.
\end{thm}
Due to the symmetry of the density association measure, we have the following corollary.

\begin{conclusion}\label{co:denexp}
If $X$ or $Z$ is binary, condition (3) in Theorem \ref{thm:6} is redundant.
\end{conclusion}
In many randomized experiments, $X$ is binary, thus we do not need to evaluate condition (3) in Theorem \ref{thm:6}.

\vspace{5mm} \noindent {\bf 5. Illustrations}

In practice, it is possible that we have two or more studies conducted
in different periods, locations or sub-populations.
When not all variables are observed in each study,
or the data may come from a conditional sampling (e.g.
from sub-populations given $Y$),
generally we can obtain only local or sub-population conclusions.
However, we are interested in global conclusions about the whole population.
We shall illustrate how to apply our theoretical results to obtain global
conclusions using two data sets.

\vspace{5mm} \noindent{\bf 5.1 An Observational Study: The National Longitudinal Surveys}

The National Longitudinal Surveys are a set of surveys designed to gather information at multiple time points on the labor market activities and other significant life events of several groups of men and women (Toomet and Henningsen, 2008).
Define $X = 1$ if a subject graduated from college, and $X=0$ otherwise.
Define $Y = 1$  if a subject belonged to a union,
and $Y = 0$ otherwise. Let $Z$ be the log of the wage.

To illustrate our theoretical results,
suppose that the data set is collected from two groups,
one is the union group with $Y=1$ and the other is from
the non-union group with $Y=0$.
That is, the data set is sampled from the conditional
distributions $f(x,z|y)$ for $y=0$ and 1.

From the data set, we can calculate
 the differences of sample means of the log wage $Z$
between different college levels  given  $Y$:
 mean$(Z|X=1, Y=1) -$mean$(Z|X=0, Y=1)=0.4351$
($p<0.001$)
and mean$(Z|X=1, Y=0) -$mean$(Z|X=0, Y=0)=0.3328$ ($p<0.001$).
Thus the college indicator $X$ has a positive expectation association
on the log wage $Z$ for both of the groups with $Y=0$ and 1.
Since there is no information on the distribution of $Y$
from the conditional sampling given $Y$,
we cannot determine whether the college indicator $X$ also has a positive
association on the log wage $Z$ marginally.
To do that, we check the conditions in Theorem \ref{thm:8} as follows.
For condition (1),
we have the difference
of observed frequencies: $\hat{P}(X=1|Y=1)-\hat{P}(X=1|Y=0)=0.0626$ ($p<0.001$), thus we have $\partial F(y|x)/ x \leq 0$ since $X$ and $Y$ are binary;
and for condition (2),
we have the differences of sample means: mean$(Z|X=1, Y=1) -$mean$(Z|X=1, Y=0)=0.1196$ ($p<0.001$) and mean$(Z| X=0,Y=1) -$mean$(Z|X=0, Y=0)=0.2218$ ($p<0.001$).
Thus we can conclude that
there is a positive expectation
association of the college indicator $X$ on the log wage $Z$ in the population.

Using the data set
containing the joint distribution of $(X,Y,Z)$,
 we obtain mean$(Z|X=1) -$mean$(Z|X=0)=0.4211$, thus confirming our conclusion.
 If we would like to assume that $Z$ follows an exponential family distribution
 conditional on $X$, then we can deduce from Corollary 2 that the density and distribution association between $X$ and $Z$ are positive.

\vspace{5mm} \noindent {\bf 5.2 A Randomized Experiment: The Job Search Intervention Study}

The Job Search Intervention Study
was a randomized field experiment that investigated
the efficacy of a job training intervention on unemployed workers
(Vinocur and Schul, 1997; Tingley et al., 2012).
The program was designed not only to increase reemployment among the unemployed, but also to enhance the mental health of the job seekers. In the JOBS II field experiment, 1,801 unemployed workers were randomly assigned to a treatment group ($X=1$) and a control group ($X=0$). Those in the treatment group participated in job-skills workshops. In the workshops, respondents learned job-searching skills and coping strategies for dealing with setbacks in the job-search process. Those in the control group received a booklet describing job-searching tips. In the follow-up interviews, two key outcome variables were measured: a continuous variable $Y$ denoting the job-search self-efficacy,
and a continuous variable $Z$ denoting depressive symptoms
based on the Hopkins Symptom Checklist.

We randomly split the data set into two subsets with the same number of observations,
and suppose that only $X$ and $Y$ were observed in the first subset
and only $Y$ and $Z$ are observed in the second subset.
The goal is to qualitatively evaluate the effect of the treatment $X$ on the depressive symptoms $Z$.
We assume the linear regression model: $E(Z'|x,y) = \beta_0 + \beta_1 x + \beta_2 y$ and $\beta_1 >0$,
where $Z'=-Z$.
This was confirmed by the full data set ($\hat{\beta}_1=0.0408$ with $p$-value=0.36).

For the condition (1) of Corollary \ref{co4},
the linear model $E(Y|x) = \beta_3 + \beta_4 x$ is saturated
because $X$ is binary,
and we obtain $\hat{\beta}_4=0.1416$ with $p$-value=0.023 from the first subset containing
$X$ and $Y$.

\begin{figure}[!h]
\centering
\includegraphics[width = \textwidth]{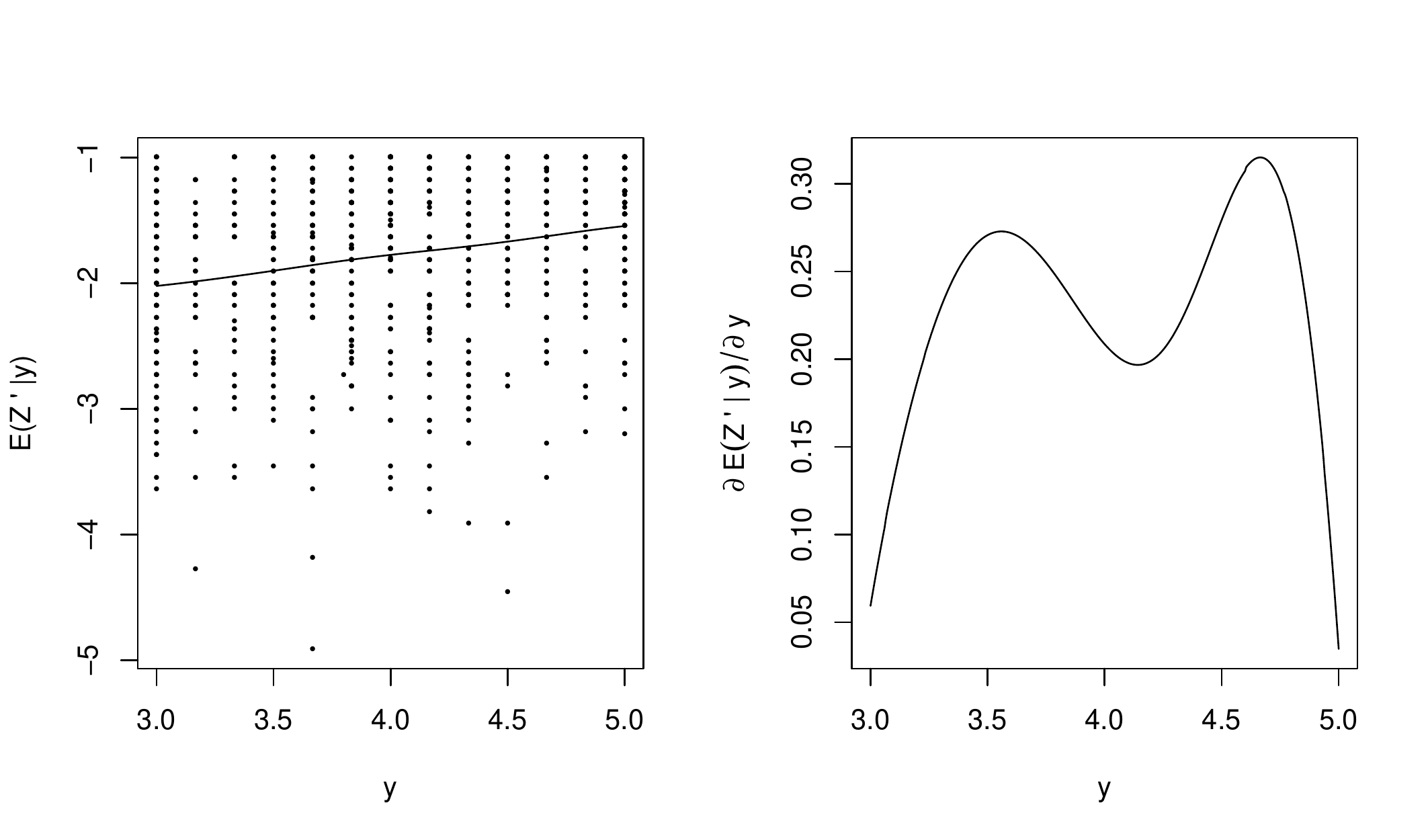}
\caption{The curves of $E(Z'|y)$ and $\partial E(Z'|y)/ \partial y$
estimated by local polynomial regression with Gaussian kernel and bandwidth $0.4$.}
\label{fg::deriv}
\end{figure}

To check condition (2)
 of Corollary \ref{co4}, we use both parametric and nonparametric approaches.
First, we use the linear model $E(Z'|y)=\beta_5 + \beta_6y$
and obtain the estimate of $\partial E(Z'|y)/ \partial y =\beta_6$:
$\hat{\beta}_6=0.2882$ ($p$-value $<$0.001)
from the second data subset of observed variables $Y$ and $Z'$.
 Next we use local polynomial regression to get nonparametric estimations of both $E(Z'|y)$ and $\partial E(Z'|y)/ \partial y$, which are shown in Figure \ref{fg::deriv}.
The scatterplot of $Y$ and $Z'$ is also shown in the left side of Figure \ref{fg::deriv}.
It is seen that $\partial E(Z'|y)/ \partial y \geq 0$ for any $y$. Both approaches confirm condition (2) of Corollary \ref{co4}, thus we conclude that $X$ has a non-negative association with $Z'$.

 From the full data set, we obtain that $E(Z|x=1) - E(Z|x=0)=0.06335$  ($p$-value $=$0.17), thus confirming our conclusion.
Although the $p$-value for $\beta_1$ is not significant, this example is useful for illustrating Corollary \ref{co4}, because  Corollary \ref{co4}  is about only the signs of parameters and does not require the significance of the parameters. We also learn from this example that the transitivity of hypothesis testing is much more difficult than the transitivity of association signs.

\vspace{5mm}
\noindent {\bf 6. Discussion and extensions}

We have discussed the transitivity of the signs of four association measures
with or without the conditional independence assumption.
These four association measures have different stringencies,
and more stringent association measures have stronger transitivity.
Consequently, the signs or directions of stringent association measures
are more easily kept for transitivity under the conditional independence assumption.
We proposed conditions for the transitivity of association measures,
and showed that these conditions are necessary and sufficient
if the intermediate variable $Y$ is from an exponential family  given  $X$.
The conditions can be checked by data based on marginal distributions
or conditional distributions from different studies.

%


The problem of transportability (Frangakis and Rubin, 2002; Rubin, 2004; Pearl and Bareinboim, 2011)
is related to transitivity discussed in this paper.
Transportability refers to two populations
sharing some features in common,
but transitivity refers to variables in a single population
which are not observed jointly.
If we treat these variables as from different populations,
then transitivity can be viewed as transportability.
Pearl and Bareinboim's approaches for transportability are for quantitative inference,
and they require more restrictive conditions
for the similarity between different populations, which may not be satisfied in practice. Qualitative reasoning is often used in researche studies and daily life, and it does not require such restrictive conditions.

 As pointed out by a referee, several applications and generalizations beyond our paper are possible. The density association is used mostly  for binary variables, and, among continuous distribution, it has been defined explicitly only for the multivariate normal distribution. Therefore, it may be of interest
 to find other continuous distributions which have nice properties  for  the density association. In this paper, the results do not require a particular DAG, and we discuss two DAGs related to our results. In fact, there are four different types of decomposable chain graphs discussed by Drton (2009), and it is worthwhile to discuss applications to these chain graphs. 
 It may also be interesting to
  simplify these conditions for transitivity 
  using the distributional results by Roverato (2013) and the conditions for traceability of paths in regression graphs by Wermuth (2012).

\vspace{5mm} \noindent {\large\bf Acknowledgment}

The authors thank a reviewer for valuable comments. This research was supported by NSFC (11171365, 11021463, 10931002).

\vspace{5mm}
\noindent {\large \bf Appendix}

\noindent {\bf Lemma 1.}
If  $h(y,a,r)$ is non-decreasing in $y$ and in $a$, and $S(y | a, r)=P(Y>y|A=a,R=r)$ is non-decreasing in $a$ for all $y$, then $E\{h(Y,A,R)|A=a,R=r\}$ is non-decreasing in $a$.
\\
 {\it Proof.}
A proof is given by VanderWeele and Robins  (2009, page 710, line 7).
As suggested by a reviewer,  we give a proof for discrete variables. Suppose $a \geq a'$ and $Y$ is discrete, taking values   $-\infty=y_0<y_1<y_2<\ldots<y_k$, then we have
\begin{eqnarray*}
&&E\{h(Y,A,R)| A=a, R=r\}-E\{h(Y,A,R)| A=a', R=r\}\\
&=& \sum_{i=1}^{k} h(y_i,a,r) P(Y=y_i | A=a, R=r)-\sum_{i=1}^{k} h(y_i,a',r) P(Y=y_i | A=a', R=r)\\
&=& \sum_{i=1}^{k} h(y_i,a,r) \{ S(y_{i-1} | a, r)- S(y_i | a, r)\}-\sum_{i=1}^{k} h(y_i,a',r) \{ S(y_{i-1} | a', r)- S(y_i | a', r)\}\\
&=& \sum_{i=1}^{k} h(y_i,a,r) \{ S(y_{i-1} | a, r)- S(y_{i-1} | a', r)\}+\sum_{i=1}^{k} \{h(y_i,a,r)- h(y_i,a',r)\} S(y_{i-1} | a', r)\\
&&- \sum_{i=1}^{k} h(y_i,a,r) \{ S(y_i | a, r)- S(y_i | a', r)\}- \sum_{i=1}^{k} \{h(y_i,a,r)- h(y_i,a',r)\} S(y_i | a', r)\\
&=&\sum_{i=2}^{k} \{h(y_i,a,r)-h(y_{i-1},a,r)\} \{ S(y_{i-1} | a, r)- S(y_{i-1} | a', r)\}\\
&&+ \sum_{i=1}^{k} \{h(y_i,a,r)- h(y_i,a',r)\} \{S(y_{i-1} | a', r)-S(y_i | a', r)\}.
\end{eqnarray*}
The final expression is non-negative since all differences in brackets are non-negative for $a \geq a'$.\\

\noindent     {\it Proof of Theorem \ref{th1}.}
We need only to prove that
$\partial \ln f(z|x)/ \partial x \geq \partial \ln f(z'|x)/\partial x$
for all $z>z'$.
When $X$ is continuous, we deduce from $ X \ind Z|Y$ that
\begin{eqnarray} \label{(A.1)}
\quad \frac{\partial \ln f(z|x)}{\partial x}
    &&=\frac{\partial  f(z|x)}{\partial x}/f(z|x)
    = \frac{\partial}{\partial x} \left\{ \int_{-\infty}^{+\infty}f(z, y|x)dy \right\} /
    f(z|x)  \nonumber \\
    &&= \int_{-\infty}^{+\infty}\frac{\partial f(y|x)}{\partial x} \frac{f(z|y)}{f(z|x)}dy
    =\int_{-\infty}^{+\infty}\frac{\partial \ln f(y|x)}{\partial x} \frac{f(z|y)f(y|x)}{f(z|x) }  dy \nonumber \\
    &&=\int_{-\infty}^{+\infty}\frac{\partial \ln f(y|x)}{\partial x}  f(y|x,z) dy  \nonumber \\
    &&=
    E\left\{ \frac{\partial \ln f(Y|x)}{\partial x}|X=x,Z=z\right\}.
\end{eqnarray}

 From $\partial^2 \ln f(y |x)/\partial y\partial x \geq 0$, we know that $\partial \ln f(y|x)/\partial x$ is non-decreasing in $y$.
Again from $ X \ind Z|Y$, we have $\ln f(x,y,z) = \ln f(y) + \ln f(z|y) + \ln f(x|y)$. By  condition (2) in Theorem \ref{th1}, we obtain
$$
    \frac{\partial^2 \ln f(y,z |x)}{\partial y\partial z}
    = \frac{\partial^2 \ln f(z|y, x)}{\partial y\partial z}
    =   \frac{\partial^2 \ln f(z|y)}{\partial y\partial z}
    = \frac{\partial^2 \ln f(y, z)}{\partial y\partial z} \geq 0.
$$
From Property 1, we get $\partial F(y|z, x)/\partial z \leq 0$, and thus
$P(Y>y|X=x,Z=z)$ is non-decreasing in $z$ for all $y$.
Applying Lemma 1 to equation (\ref{(A.1)}),
we conclude that $\partial \ln f(z|x)/\partial x$ is non-decreasing in $z$.

When $X$ is discrete, we need only to prove that, for all $z>z'$,
$$
    \frac{ f(z|x=1)}{f(z|x=0)}\geq\frac{ f(z'|x=1)}{f(z'|x=0)},
    $$
or, equivalently,
$$
    \frac{ f(z|x=1)-f(z|x=0)}{f(z|x=0)}\geq\frac{ f(z'|x=1)-f(z'|x=0)}{f(z'|x=0)}.
    $$
We compute that
\begin{eqnarray*}
 \frac{ f(z|x=1)-f(z|x=0)}{f(z|x=0)}
    &=& \int_{-\infty}^{+\infty}\frac{ f(z|y) \left\{ f(y|x=1)-f(y|x=0)  \right\} }{f(z|x=0)}dy  \\
    &=&\int_{-\infty}^{+\infty}\frac{f(y|x=1)-f(y|x=0)}{f(y|x=0)}\frac{f(z|y)f(y|x=0)}{f(z|x=0) }  dy  \\
    &=&\int_{-\infty}^{+\infty}\frac{f(y|x=1)-f(y|x=0)}{f(y|x=0)} f(y|z, x=0)dy\\
    &=&E\left\{ \frac{f(Y|x=1)-f(Y|x=0)}{f(Y|x=0)}|X=0,Z=z \right\}.
\end{eqnarray*}
From  conditions (1) and (2) in Theorem \ref{th1},
we have that $\{f(y|x=1)-f(y|x=0)\}/f(y|x=0) $ is non-decreasing in $y$ and that $P(Y>y|X=x,Z=z)$ is non-decreasing in $z$ for all $y$.
Thus by Lemma 1, we conclude that $f(z|x=1) / f(z|x=0)  \geq  f(z'|x=1) / f(z'|x=0)$.\\

\noindent {\it Proof of Theorem \ref{th2}.}
By $X \ind Z|Y$, we have
\begin{eqnarray*}
  F(z|x)
= \int_{-\infty}^{+\infty}   F(z| y) F(dy|x)= E\left\{   F(z|Y)|X=x \right\}.
\end{eqnarray*}
From  conditions (1) and (2) in Theorem \ref{th2}, $F(z|y)$ is non-increasing in $y$, and $P(Y>y|X=x)$ is non-decreasing in $x$ for all $y$. By Lemma 1, $F(z|x)$ is non-increasing in $x$.\\


\noindent {\it Proof of Theorem \ref{co2}.}
For the exponential family, we
have that $\partial^2 \ln f(x,y) /\partial x \partial y  =  (\partial \theta_x / \partial x) /a(\phi)$ and
$\partial E(Y|x) / \partial x = \partial b'(\theta_x) / \partial x = b''(\theta_x) (\partial \theta_x / \partial x) = \text{var}(Y|x)\times \\ (\partial \theta_x / \partial x)/a(\phi)$. Thus we obtain that $\partial^2 \ln f(x,y) / \partial x \partial y$ and $\partial E(Y|x) / \partial x$ have the same sign, which implies the conclusion. \\

\noindent    {\it Proof of Corollary \ref{co3}.}
The implication relationships from the signs of association measures between $X$ and $Y$ to the signs of association measures between $X$ and $Z$ can be deduced from Theorems \ref{th1} to \ref{co2}.
Below we show three implication relationships
from the signs of measures between $X$ and $Z$ to the signs of measures between $X$ and $Y$.
From result (1) of Theorem \ref{co2}, we need to show that
$E(Z|x)$ increasing in $x$ implies
$E(Y|x)$ increasing in $x$. By $X\ind Z|Y$, we have from the proof of Theorem \ref{th3} that
\begin{eqnarray*}
E(Z|x) -E(Z|x')= - \int_{-\infty}^{+\infty}
\frac{\partial E(Z|y)}{\partial y} \{ F(y|x) - F(y|x')\}  dy.
\end{eqnarray*}
We use proof by contradiction;
suppose that there exists $x>x'$ such that $E(Y|x)<E(Y|x')$.
Then from the property of the exponential family in Theorem \ref{co2},
we have
$F(y|x)>F(y|x')$ for all $y$.
Because $\partial E(Z|y) / \partial y$ is strictly positive
for a non-zero measure set, we get
that $E(Z|x) -E(Z|x')<0$ from Theorem \ref{th3}, which contradicts the condition
of a non-negative association between $X$ and $Z$.

Results (2) and (3) of Corollary \ref{co3}
can be obtained immediately from the above result (1) and Theorem \ref{co2}.
\\

\noindent {\it Proof of Theorem \ref{thm:6}.} We need only to prove that $\partial \ln f(z|x)/\partial x$ is non-decreasing in $z$.
When $X$ is continuous, for $z>z'$, we have
    \begin{eqnarray*}
    && \quad \frac{\partial \ln f(z|x)}{\partial x} =\frac{\partial f(z|x)}{\partial x}/f(z|x)
    =\left\{\frac{\partial}{\partial x} \int_{-\infty}^{+\infty}f(z|y, x)f(y|x)dy \right\}/
    f(z|x)\\
    &&=\int_{-\infty}^{+\infty}   \left \{ \frac{\partial f(z|y, x)}{\partial x}\cdot \frac{f(y|x)}{f(z|x)}
    +\frac{\partial f(y|x)}{\partial x}\cdot \frac{f(z|y, x)}{f(z|x)}   \right\}    dy\\
    &&=\int_{-\infty}^{+\infty}  \left\{     \frac{\partial f(z|y, x)}{\partial x}/f(z|y, x)
    +\frac{\partial f(y|x)}{\partial x}/f(y|x)   \right\}   \cdot \frac{f(y|x)f(z|y, x)}{f(z|x)}dy\\
    &&=\int_{-\infty}^{+\infty} \left\{\frac{\partial \ln f(z|y, x)}{\partial x} +
    \frac{\partial \ln f (y|x)}{\partial x}\right\} f(y|x, z)dy.\\
    &&=E\left\{   \frac{\partial \ln f(z|Y, x)}{\partial x}|Z=z,X=x\right\}
    + E\left\{  \frac{\partial \ln f (Y|x)}{\partial x}|Z=z,X=x\right\}.
    \end{eqnarray*}
From  the assumption and condition (3) in Theorem \ref{thm:6}, $\partial \ln f(z|Y, x)/ \partial x$ is non-decreasing in $y$ and $z$; from condition (1) in Theorem \ref{thm:6}, $ \partial \ln f (Y|x) / \partial x$ is nondecreasing in $y$;  and from  condition (2) in Theorem \ref{thm:6},  $P(Y>y|X=x,Z=z)$ is non-decreasing in $z$ for all $y$. By Lemma 1, we have that
$$E\left\{ \frac{ \partial \ln f(z|Y, x)}{  \partial x}|Z=z,X=x\right\}+E\left\{ \frac{ \partial \ln f (Y|x)}{ \partial x}|Z=z,X=x\right\}$$
 is non-decreasing in $z$.
 
When $X$ is discrete, we need only to prove that, for $z>z'$,
$$
    \frac{ f(z|x=1)}{f(z|x=0)}\geq\frac{ f(z'|x=1)}{f(z'|x=0)}.
    $$
We have that
\begin{eqnarray*}
    && \frac{ f(z|x=1)}{f(z|x=0)}
    = \int_{-\infty}^{+\infty}\frac{ f(z|y,x=1)  f(y|x=1)}{f(z|x=0)}dy  \\
    &=&\int_{-\infty}^{+\infty}\frac{f(z|y,x=1)}{f(z|y,x=0)}\frac{f(y|x=1)}{f(y|x=0)}  \frac{f(z|y,x=0)f(y|x=0)}{f(z|x=0)}dy  \\
    &=&\int_{-\infty}^{+\infty}\frac{f(z|y,x=1)}{f(z|y,x=0)}\frac{f(y|x=1)}{f(y|x=0)} f(y|z,x=0)dy  \\
    &=&E\left\{   \frac{f(z|Y,x=1)}{f(z|Y,x=0)}\frac{f(Y|x=1)}{f(Y|x=0)}|X=0,Z=z \right\}.
\end{eqnarray*}

From the assumption and  condition (1) in Theorem \ref{thm:6},
$\{  f(z|y,x=1) f(y|x=1) \} / \{ f(z|y,x=0)  f(y|x=0)\}$
is non-decreasing in $y$ and $z$. From  condition (2) in Theorem \ref{thm:6}, $P(Y>y|X=0,Z=z)$ is non-decreasing in $z$ for all $y$.
Therefore, we have $f(z|x=1) / f(z|x=0)  \geq   f(z'|x=1) / f(z'|x=0)$.
\\

\noindent   {\it Proof of Theorem \ref{thm:7}.}
For $F(z|x)=E\left\{ F(z|Y,x)|X=x\right\} $, we have from the assumption and  condition (2) in Theorem \ref{thm:7} that $F(z|y,x)$ is non-increasing in $y$ and $x$. From  condition (1) in Theorem \ref{thm:7}, $P(Y>y|X=x)$ is non-decreasing in $x$ for all $y$. Thus we have that  $\partial F(z|x) / \partial x \leq 0$.
\\

\noindent  {\it Proof of Theorem \ref{thm:8}.}
Since $E(Z|x)=E\left\{ E(Z|Y,x)|X=x \right\}$, we prove that $\partial E(z|x) / \partial x \leq 0$ using a similar argument as in the proof of Theorem 5.
\\

\noindent    {\it Proof of Theorem \ref{thm:exp}.}
From Theorem \ref{thm:8}, we have $\partial E(Z| x) / \partial x \geq 0, \forall x$, and then from Theorem \ref{co2}, we have $\partial^2 \ln f(x, z) / \partial x\partial z \geq 0$, $\forall x, z$.\\

\noindent  {\it Proof of Corollary \ref{co:simpson}.}
We  prove this only for Theorem \ref{thm:6}.  Obviously, the assumption $\partial^2 \ln f(x, z|y) / \partial x\partial z \geq 0$ can be evaluated by $f(x,z | y)$. Condition (1) in Theorem \ref{thm:6} can be evaluated by $f(x | y)$, which can be obtained after marginalizing $f(x,z | y)$ over $z$. For conditions (2) and (3), we can rewrite them as $\partial^2 \ln f( z,x | y)/\partial y\partial z\geq 0$ and $\partial^2 \ln f(x,z | y)/\partial x \partial y\geq 0$ respectively. Therefore, the assumption and conditions can all be evaluated by $f(x,z | y)$.\\

\noindent  {\it Proof of Corollary \ref{co4}.}
From the linear model, we have $E(Z|x)=\beta_{0}+\beta_{1}x+\beta_{2}E(Y|x)$ and
$\partial E(Z|x) / \partial x  = \beta_1 + \beta_2 \partial E(Y|x)/\partial x   = \beta_{1}+\beta_{2}\beta_4$.
Thus, we have $\partial E(Z|x) / \partial x \geq 0$
if $\beta_1$, $\beta_2$ and $\beta_4$ are non-negative.

Suppose $a=-\beta_{2} / \beta_{1}$, and we need only to prove the result
for the case that $\beta_2<0$ but $\partial E(Z|Y=y) / \partial y \geq 0$.
We use proof by contradiction, and
suppose that $\partial E(Z|x) / \partial x < 0$ for some $x$.
Then we have
$\partial E(Y|x) / \partial x = \beta_{4} > -\beta_{1} / \beta_{2} =1 /a$.
Since $\partial E(Z|y) / \partial y=\beta_{1} \partial E(X|y) / \partial x +\beta_{2} \geq 0$,
we have $\partial E(X|Y=y) / \partial y \geq -\beta_{2} / \beta_{1} = a$.
From the linear model of $Y$, we have $\partial E(Y-\beta_{4}X|X=x) / \partial x =0$.
We deduce that
\begin{equation} \label{(i)}
\text{cov}(X,Y) = \beta_{4} \text{var}(X) > \text{var}(X)/a.
\end{equation}

Define $b=\inf_y\{  \partial E(X|Y=y) / \partial y\} $.
We have $b \geq a$ and $\partial E(X-bY|Y=y) / \partial y \geq b - b = 0$.
From Property 1, we get $\text{cov}(X-bY, Y)\geq 0$.
Thus we obtain
\begin{equation} \label{(ii)}
\text{cov}(X,Y)= \text{cov}(X-bY,Y) + b \text{var}(Y) \geq b \text{var}(Y) \geq a \text{var}(Y).
\end{equation}
 From equations (\ref{(i)}) and (\ref{(ii)}),
we have $\text{cov}(X,Y) > \left\{   \text{var}(X)\text{var}(Y) \right\}^{1/2}$, which is impossible
since the correlation coefficient cannot be larger than $1$.\\

\noindent  {\it Proof of Corollary \ref{co:ind}.}
We first prove results (2) and (3).
Since $F(z|x)=E\left\{ F(z|Y,x)|X=x\right\}= E_{Y}\left\{ F(z|Y,x)\right\}$, and $E(Z|x)=E\left\{ E(Z|Y,x) | X=x \right\}=E_Y\left\{ E(Z|Y,x)  \right\}$, we only need $\partial F(z| y,x)/\partial x \leq 0, \forall x,y,z$ for Theorem \ref{thm:7} and $\partial E(Z| y,x)/\partial x \geq 0, \forall x,y$ for Theorem \ref{thm:8}. For result (1), according to Theorem \ref{co2}, when $X$ or $Z$ is binary, the density association is equivalent to the expectation association, thus we need only $\partial^2 f(x,z|y)/\partial x \partial z \geq 0, \forall x,y,z$.


\vspace{5mm}
\noindent{\large\bf References}

\begin{description}

\item
\textsc{Appleton, D. R., French, J. M. } and \textsc{Vanderpump, M. P. J.} (1992). Ignoring a covariate: an example of Simpson's paradox.
\textit{Am. Stat.}
\textbf{50}, 340--341.

\item
\textsc{Birch, M. W.} (1963). Maximum likelihood in three-way contingency tables.
\textit{J. Roy. Statist. Soc. Ser. B}
\textbf{25}, 220--233.

\item
\textsc{Chen, H., Geng, Z.} and \textsc{Jia, J.} (2007). Criteria for surrogate end points.
\textit{J. Roy. Statist. Soc. Ser. B}
\textbf{69}, 919--932.

\item
\textsc{Cochran, W. G. } (1938). The omission or addition of an independent variate in multiple linear regression.
\textit{Supp. J. Roy. Statist. Soc.}
\textbf{5}, 171--176.

\item
\textsc{Cox, D. R.} and \textsc{Wermuth, N.}  (2003). A general condition for avoiding effect reversal after marginalization.
\textit{J. Roy. Statist. Soc. Ser. B}
\textbf{48}, 197--205.

\item
\textsc{Drton, M. } (2009). Discrete chain graph models.
\textit{Bernoulli}
\textbf{15}, 736--753.

%

\item
\textsc{Frangakis, C. E.} and \textsc{Rubin, D. B.}  (2002).
Principal stratification in causal inference.
\textit{Biometrics}
\textbf{58}, 21--29.

%

\item
\textsc{Ju, C.} and \textsc{Geng, Z.} (2010).
Criteria for surrogate end points based on causal distributions. \textit{J. Roy. Statist. Soc. Ser. B} \textbf{72}, 129--142.

%
%

\item
\textsc{Ln\v{e}ni\v{c}ka, R.} and \textsc{Mat\'{u}\v{s}, F.} (2007).
On Gaussian conditional independence structures. \textit{Kybernetika} \textbf{43}, 323--342.

\item
 \textsc{Pearl, J.} and \textsc{Bareinboim, E.} (2011).
Transportability of causal and statistical relations: A formal approach. {\it Data Mining Workshops (ICDMW), 2011 IEEE 11th International Conference on} IEEE, 540--547.

\item
\textsc{Prentice, R. L.} (1989).
Surrogate endpoints in clinical trials: definition and operational criteria.
 \textit{Stat. Med.} \textbf{8}, 431--440.


\item
\textsc{Roverato, A.}  (2013).
Dichotomization invariant log-mean linear parameterization for discrete graphical models of marginal independence. \href{http://arxiv.org/pdf/1302.4641.pdf}{arXiv:1302.4641}.

\item
\textsc{Rubin, D. B. } (2004).
Direct and indirect causal effects via potential outcomes.
\textit{Scand. J. Stat.} \textbf{31}, 161--170.

\item
\textsc{Toomet, O.} and \textsc{Henningsen, A.} (2008).
Sample selection models in R: package
  sample Selection.
\textit{J. Stat. Softw.}
\textbf{27},  \\
\url{http://www.jstatsoft.org/v27/i07/}.

\item
\textsc{Tingley, D., Yamamoto, T., Keele, L.} and \textsc{Imai, K.} (2012).
mediation: R package for causal mediation analysis.
R package version 4.2.
\\  \url{http://CRAN.R-project.org/package=mediation} .

\item
\textsc{Vanderweele, T. J.} and \textsc{Robins, J. M.} (2009).
Properties of monotonic effects on
directed acyclic graphs.
\textit{J. Mach. Learn. Res.}
\textbf{10}, 699--718.

\item
\textsc{Vanderweele, T. J.} and \textsc{Robins, J. M.} (2010).
Signed directed acyclic graphs for causal inference.
\textit{J. Roy. Statist. Soc. Ser. B}
\textbf{72}, 111--127.

\item
\textsc{Vanderweele, T. J.} and \textsc{Tan, Z.} (2012).
Directed acyclic graphs with edge-specific bounds.
\textit{Biometrika}
\textbf{99}, 115--126.

\item
\textsc{Vinokur, A. } and \textsc{Schul, Y.} (1997).
Mastery and inoculation against setbacks as active ingredients in the jobs intervention for the unempllyed. \textit{J. Consult. Clin. Psych.} \textbf{65}, 867--877.

\item
\textsc{Wermuth, N. }  (2012).
Traceable regressions. \textit{Int. Stat. Rev.} \textbf{80}, 415-438.

\item
\textsc{Whittaker, J.} (1990).
\textit{Graphical Models in Applied
Multivariate Statistics.} New York: John Wiley \& Sons.

%
%
%

\item
\textsc{Xie, X., Ma, Z.} and \textsc{Geng, Z.} (2008).
Some association measures and their collapsibility. \textit{Statistica Sinica} \textbf{18}, 1165--1183.

\end{description}


\vskip .65cm
\noindent
School of Mathematical Sciences,
Peking University, Beijing 100871, China
\vskip 2pt
\noindent
E-mail: (zhichaojiang@pku.edu.cn)
\vskip 2pt

\noindent
Department of Statistics, Harvard University,
Cambridge 02138, U.S.A.
\vskip 2pt
\noindent
E-mail: (pengding@fas.harvard.edu)
\vskip 2pt

\noindent
Center for Statistical Science,
School of Mathematical Sciences,
Peking University, Beijing 100871, China
\vskip 2pt
\noindent
E-mail: (zhigeng@pku.edu.cn)
\vskip .3cm


\end{document}